\documentclass{amsart}
\usepackage{amssymb}

\usepackage{amscd}
\usepackage{thmdefs}

\input{tcilatex}
\theoremstyle{definition}
\theoremstyle{remark}
\numberwithin{equation}{section}

\input tcilatex

\begin{document}
\title[Graph $W^{*}$-Probability on $L(F_{N})$ ]{Graph $W^{*}$-Probability on the Free Group Factor $L(F_{N})$}
\author{Ilwoo Cho}
\address{Dep. of Math, Univ. of Iowa, Iowa City, IA, U. S. A}
\email{ilcho@math.uiowa.edu}
\maketitle

\begin{abstract}
In this paper, we will consider the free probability on the\ free group
factor $L(F_{N}),$ in terms of Graph $W^{*}$-probability, where $F_{k}$ is
the free group with $k$-generators. The main result of this paper is to
reformulate the moment series and the R-transform of the operator $%
T_{0}=\sum_{j=1}^{N}\left( x_{j}\right) ,$ where $x_{j}$'s are free
semicircular elements, $j=1,...,N,$ generating $L(F_{N}),$ by using the
Graph $W^{*}$-probability technique. To do that, we will use the graph $%
W^{*} $-probability technique to find the moments and cumulants of the
identically distributed random variable $T$ with $T_{0}.$ This will be a
good example of an application of Graph $W^{*}$-Probability Theory. We also
can see how we can construct the $W^{*}$-subalgebra which is isomorphic to
the free group factor $L(F_{N})$ in the graph $W^{*}$-probability space $%
\left( W^{*}(G),E\right) .$
\end{abstract}

\strut

\strut In this paper, we will reformulate the moments and cumulants of the
so-called generating operator $T_{0}$ $=$ $\sum_{j=1}^{N}x_{j}$, where $%
x_{1},$ ..., $x_{N}$ are free semicircular elements generating the free
group factor $L(F_{N}),$ in terms of the graph $W^{*}$-probability theory
considered in [14], [15], [16] and [17]. Voiculescu showed that the free
group factor $L(F_{N})$ is generated by $N$-semicircular elements which are
free from each other (See [9]). The moments and cumulants of such elements
are known but we will recompute them, by using the graph $W^{*}$-probability
technique. We will construct a $D_{G}$-semicircular subalgebra $S_{G}$ which
is isomorphic to the free group factor $L(F_{N})$ embedded in a certain
graph $W^{*}$-algebra $W^{*}(G).$ And, by constructing an operator $T$ which
is identically distributed with $T_{0}$ in $S_{G}$, we will recompute the
moments and cumulants of $T_{0}.$ This would be the an application of graph $%
W^{*}$-probability theory. Also, we will embed the free group factor $%
L(F_{N})$ into an arbitrary graph $W^{*}$-probability space $\left(
W^{*}(G),E\right) ,$ where $G$ is a countable directed graph containing at
least one vertex $v_{0}$ having $N$-basic loops concentrated on $v_{0}.$

\strut

Let $G$ be a countable directed graph and $W^{*}(G),$ the graph $W^{*}$%
-algebra. By defining the diagonal subalgebra $D_{G}$ and the canonical
conditional expectation $E$ $:$ $W^{*}(G)$ $\rightarrow $ $D_{G},$ we can
construct the graph $W^{*}$-probability space $\left( W^{*}(G),\text{ }%
E\right) $ over $D_{G},$ as a $W^{*}$-probability space with amalgamation
over $D_{G}$. All elements in $(W^{*}(G),$ $E)$ are called $D_{G}$-valued
random variables. The $D_{G}$-freeness is observed in [14] and [15]. The
generators $L_{w_{1}}$ and $L_{w_{2}}$ are free over $D_{G}$ if and only if $%
w_{1}$ and $w_{2}$ are diagram-distinct, in the sense that they have
different diagram on the graph $G,$ graphically. There are plenty of
interesting examples of $D_{G}$-valued random variables in this structure
(See [15]), including $D_{G}$-semicircular elements, $D_{G}$-valued
R-diagonal elements and $D_{G}$-even elements.

\strut

In this paper, we will regard the free group factor $L(F_{N})$ as an
embedded $W^{*}$-subalgebra of the graph $W^{*}$-algebra $W^{*}(G)$, where $%
G $ is a directed graph with

\strut

\begin{center}
$V(G)=\{v\}$ \ \ and \ \ $E(G)=\{l_{1},...,l_{N}\},$
\end{center}

\strut

where $l_{j}=vl_{j}v$ is a loop-edge, for all $j=1,...,N.$ Then the
generating operator $T_{0}$ of $L(F_{N})$ is identically distributed with
the operator $T$ in $W^{*}(G)$ such that

\strut

\begin{center}
$T=\sum_{j=1}^{N}\frac{1}{\sqrt{2}}\left( L_{l_{j}}+L_{l_{j}}^{*}\right) .$
\end{center}

\strut

Notice that, since $l_{j}$'s are diagram-distinct, in the sense that they
have mutually different diagrams in the graph $G,$ $L_{l_{j}}+L_{l_{j}}^{*}$%
's are free from each other over $D_{G}=\Bbb{C}$ in $\left(
W^{*}(G),E\right) ,$ for \ $j=1,...,N.$ Futhermore, by [14], we know that
the summands $L_{l_{j}}+L_{l_{j}}^{*}$'s are all $D_{G}=\Bbb{C}$%
-semicircular. So, we can get the cumulants of $T$ somewhat easily. And the
moments and cumulants gotten from this are same as those of the generating
operator $T_{0}$ in the free group factor $L(F_{N}).$

\strut

Recall that studying moment series of a random variable is studying
distribution of the random variable. By the moment series of the random
variable, we can get the algebraic and combinatorial information about the
distribution of the random variable. Also, alternatively, the R-transforms
of random variables contains algebraic and combinatorial information about
distributions of the random variables. So, to study moment series and
R-transforms of random variables is very important to study distributions of
those random variables. Moreover, to study R-transform theory allows us to
understand the freeness of random variables. This paper deals with the
operator $T_{0}$ $=$ $\sum_{j=1}^{N}$ $x_{j}$ of the free group factor $%
L(F_{N})$, where $x_{1},$ $...,$ $x_{N}$ are semicircular elements, which
are free from each other, generating $L(F_{N})$. Notice that, by Voiculescu,
we can regard the free group factor $L(F_{N})$ as the von Neumann algebra $%
vN $ $\left( \{x_{j}\}_{j=1}^{N}\right) .$ We will consider this free
probabilistic data about $T_{0},$ by using the graph $W^{*}$-probability
technique.

\strut

In Chapter 1, we will review the graph $W^{*}$-probability theory. In
Chapter 2, we will construct the graph $W^{*}$-probability space induced by
the one-vertex-$N$-loop-edge graph and define a $W^{*}$-subalgebra of the
graph $W^{*}$-algebra, which is isomorphic to $L(F_{N}),$ in the sense of
[9]. In Chapter 3, we will re-compute the moments and cumulants of the
generating operator $T_{0}$ of $L(F_{N}),$ by using the Graph $W^{*}$%
-probability technique. In Chapter 4, we will observe the embeddings of $%
L(F_{N})$ into $W^{*}(G),$ where $G$ is an arbitrary countable directed
graph containing a vertex $v_{0}$ having $N$-basic loops concentrate on it.

\strut

\strut

\strut

\section{Graph $W^{*}$-Probability Spaces}

\strut

\strut

Let $G$ be a countable directed graph and let $\Bbb{F}^{+}(G)$ be the free
semigroupoid of $G.$ i.e., the set $\mathbb{F}^{+}(G)$ is the collection of
all vertices as units and all admissible finite paths of $G.$ Let $w$ be a
finite path with its source $s(w)=x$ and its range $r(w)=y,$ where $x,y\in
V(G).$ Then sometimes we will denote $w$ by $w=xwy$ to express the source
and the range of $w.$ We can define the graph Hilbert space $H_{G}$ by the
Hilbert space $l^{2}\left( \mathbb{F}^{+}(G)\right) $ generated by the
elements in the free semigroupoid $\mathbb{F}^{+}(G).$ i.e., this Hilbert
space has its Hilbert basis $\mathcal{B}=\{\xi _{w}:w\in \mathbb{F}%
^{+}(G)\}. $ Suppose that $w=e_{1}...e_{k}\in FP(G)$ is a finite path with $%
e_{1},...,e_{k}\in E(G).$ Then we can regard $\xi _{w}$ as $\xi
_{e_{1}}\otimes ...\otimes \xi _{e_{k}}.$ So, in [10], Kribs and Power
called this graph Hilbert space the generalized Fock space. Throughout this
paper, we will call $H_{G}$ the graph Hilbert space to emphasize that this
Hilbert space is induced by the graph.

\strut

Define the creation operator $L_{w},$ for $w\in \mathbb{F}^{+}(G),$ by the
multiplication operator by $\xi _{w}$ on $H_{G}.$ Then the creation operator 
$L$ on $H_{G}$ satisfies that

\strut

(i) \ $L_{w}=L_{xwy}=L_{x}L_{w}L_{y},$ for $w=xwy$ with $x,y\in V(G).$

\strut

(ii) $L_{w_{1}}L_{w_{2}}=\left\{ 
\begin{array}{lll}
L_{w_{1}w_{2}} &  & \text{if }w_{1}w_{2}\in \mathbb{F}^{+}(G) \\ 
&  &  \\ 
0 &  & \text{if }w_{1}w_{2}\notin \mathbb{F}^{+}(G),
\end{array}
\right. $

\strut

\ \ \ \ for all $w_{1},w_{2}\in \mathbb{F}^{+}(G).$

\strut

Now, define the annihilation operator $L_{w}^{*},$ for $w\in \mathbb{F}%
^{+}(G)$ by

\strut

\begin{center}
$L_{w}^{\ast }\xi _{w^{\prime }}\overset{def}{=}\left\{ 
\begin{array}{lll}
\xi _{h} &  & \text{if }w^{\prime }=wh\in \mathbb{F}^{+}(G)\xi \\ 
&  &  \\ 
0 &  & \text{otherwise.}
\end{array}
\right. $
\end{center}

\strut

The above definition is gotten by the following observation ;

\strut

\begin{center}
$
\begin{array}{ll}
<L_{w}\xi _{h},\xi _{wh}>\, & =\,<\xi _{wh},\xi _{wh}>\, \\ 
& =\,1=\,<\xi _{h},\xi _{h}> \\ 
& =\,<\xi _{h},L_{w}^{*}\xi _{wh}>,
\end{array}
\,$
\end{center}

\strut

where $<,>$ is the inner product on the graph Hilbert space $H_{G}.$ Of
course, in the above formula we need the admissibility of $w$ and $h$ in $%
\mathbb{F}^{+}(G).$ However, even though $w$ and $h$ are not admissible
(i.e., $wh\notin \mathbb{F}^{+}(G)$), by the definition of $L_{w}^{\ast },$
we have that

\strut

\begin{center}
$
\begin{array}{ll}
<L_{w}\xi _{h},\xi _{h}> & =\,<0,\xi _{h}> \\ 
& =0=\,<\xi _{h},0> \\ 
& =\,<\xi _{h},L_{w}^{*}\xi _{h}>.
\end{array}
\,\,$
\end{center}

\strut

Notice that the creation operator $L$ and the annihilation operator $L^{*}$
satisfy that

\strut

(1.1) \ \ \ $L_{w}^{*}L_{w}=L_{y}$ \ \ and \ \ $L_{w}L_{w}^{*}=L_{x},$ \ for
all \ $w=xwy\in \mathbb{F}^{+}(G),$

\strut

where $x,y\in V(G).$ Remark that if we consider the von Neumann algebra $%
W^{*}(\{L_{w}\})$ generated by $L_{w}$ and $L_{w}^{*}$ in $B(H_{G}),$ then
the projections $L_{y}$ and $L_{x}$ are Murray-von Neumann equivalent,
because there exists a partial isometry $L_{w}$ satisfying the relation
(1.1). Indeed, if $w=xwy$ in $\mathbb{F}^{+}(G),$ with $x,y\in V(G),$ then
under the weak topology we have that

\strut

(1,2) \ \ \ $L_{w}L_{w}^{*}L_{w}=L_{w}$ \ \ and \ \ $%
L_{w}^{*}L_{w}L_{w}^{*}=L_{w}^{*}.$

\strut

So, the creation operator $L_{w}$ is a partial isometry in $W^{*}(\{L_{w}\})$
in $B(H_{G}).$ Assume now that $v\in V(G).$ Then we can regard $v$ as $%
v=vvv. $ So,

\strut

(1.3) $\ \ \ \ \ \ \ \ \ L_{v}^{*}L_{v}=L_{v}=L_{v}L_{v}^{*}=L_{v}^{*}.$

\strut

This relation shows that $L_{v}$ is a projection in $B(H_{G})$ for all $v\in
V(G).$

\strut

Define the \textbf{graph }$W^{*}$\textbf{-algebra} $W^{*}(G)$ by

\strut

\begin{center}
$W^{*}(G)\overset{def}{=}\overline{%
\mathbb{C}[\{L_{w},L_{w}^{*}:w\in
\mathbb{F}^{+}(G)\}]}^{w}.$
\end{center}

\strut

Then all generators are either partial isometries or projections, by (1.2)
and (1.3). So, this graph $W^{\ast }$-algebra contains a rich structure, as
a von Neumann algebra. (This construction can be the generalization of that
of group von Neumann algebra.) Naturally, we can define a von Neumann
subalgebra $D_{G}\subset W^{\ast }(G)$ generated by all projections $L_{v},$ 
$v\in V(G).$ i.e.

\strut

\begin{center}
$D_{G}\overset{def}{=}W^{*}\left( \{L_{v}:v\in V(G)\}\right) .$
\end{center}

\strut

We call this subalgebra the \textbf{diagonal subalgebra} of $W^{*}(G).$
Notice that $D_{G}=\Delta _{\left| G\right| }\subset M_{\left| G\right| }(%
\mathbb{C}),$ where $\Delta _{\left| G\right| }$ is the subalgebra of $%
M_{\left| G\right| }(\mathbb{C})$ generated by all diagonal matrices. Also,
notice that $1_{D_{G}}=\underset{v\in V(G)}{\sum }L_{v}=1_{W^{*}(G)}.$

\strut

If $a\in W^{*}(G)$ is an operator, then it has the following decomposition
which is called the Fourier expansion of $a$ ;

\strut

(1.4) $\ \ \ \ \ \ \ \ \ \ \ a=\underset{w\in \mathbb{F}^{+}(G:a),\,u_{w}\in
\{1,*\}}{\sum }p_{w}^{(u_{w})}L_{w}^{u_{w}},$

\strut

where $p_{w}^{(u_{w})}\in \Bbb{C}$, $u_{w}\in \{1,*\},$ and $\mathbb{F}%
^{+}(G:a)$ is the support of $a$ defined by

\strut

\begin{center}
$\mathbb{F}^{+}(G:a)=\{w\in \mathbb{F}^{+}(G):p_{w}^{(u_{w})}\neq 0\}.$
\end{center}

\strut

Remark that the free semigroupoid $\mathbb{F}^{+}(G)$ has its partition $%
\{V(G),$ $FP(G)\},$ as a set. i.e.,

\strut

\begin{center}
$\mathbb{F}^{+}(G)=V(G)\cup FP(G)$ \ \ and \ \ $V(G)\cap FP(G)=\emptyset .$
\end{center}

\strut

So, the support of $a$ is also partitioned by

\strut

\begin{center}
$\mathbb{F}^{+}(G:a)=V(G:a)\cup FP(G:a),$
\end{center}

\strut where

\begin{center}
$V(G:a)\overset{def}{=}V(G)\cap \mathbb{F}^{+}(G:a)$
\end{center}

and

\begin{center}
$FP(G:a)\overset{def}{=}FP(G)\cap \mathbb{F}^{+}(G:a).$
\end{center}

\strut

So, the above Fourier expansion (1.4) of the random variable $a$ can be
re-expressed by

\strut

(1.5) $\ \ \ \ \ \ \ \ \ \ \ \ a=\underset{v\in V(G:a)}{\sum }p_{v}L_{v}+%
\underset{w\in FP(G:a),\,u_{w}\in \{1,*\}}{\sum }%
p_{w}^{(u_{w})}L_{w}^{u_{w}}.$

\strut

We can easily see that if $V(G:a)\neq \emptyset ,$ then $\underset{v\in
V(G:a)}{\sum }p_{v}L_{v}$ is contained in the diagonal subalgebra $D_{G}.$
Also, if $V(G:a)$ $=$ $\emptyset ,$ then $\underset{v\in V(G:a)}{\sum }%
p_{v}L_{v}$ $=$ $0_{D_{G}}.$ So, we can define the following canonical
conditional expectation $E$ $:$ $W^{*}(G)$ $\rightarrow $ $D_{G}$ by

\strut

(1.6) \ \ \ $\ \ \ \ \ \ \ \ \ \ \ \ \ \ \ \ \ \ \ \ E(a)\overset{def}{=}%
\underset{v\in V(G:a)}{\sum }p_{v}L_{v},$

\strut

for all $a\in W^{*}(G)$ having its Fourier expansion (1.5). Indeed, $E$ is a
well-determined conditional expectation. Moreover it is faithful, in the
sense that if $E(a^{*}a)$ $=$ $0_{D_{G}},$ then $a$ $=$ $0_{D_{G}},$ for $a$ 
$\in $ $W^{*}(G).$

\strut \strut \strut \strut

\begin{definition}
We say that the algebraic pair $\left( W^{*}(G),E\right) $ is the graph $%
W^{*}$-probability space over the diagonal subalgebra $D_{G}$.
\end{definition}

\strut

We will define the following free probability data of $D_{G}$-valued random
variables in $(W^{*}(G),$ $E).$

\strut

\begin{definition}
Let $W^{*}(G)$ be the graph $W^{*}$-algebra induced by $G$ and let $a\in
W^{*}(G).$ Define the $n$-th ($D_{G}$-valued) moment of $a$ by

\strut

$\ \ \ \ \ \ \ \ \ \ \ \ E\left( d_{1}ad_{2}a...d_{n}a\right) ,$ for all $%
n\in \mathbb{N}$,

\strut

where $d_{1},...,d_{n}\in D_{G}$. Also, define the $n$-th ($D_{G}$-valued)
cumulant of $a$ by

\strut

$\ \ \ \ \ k_{n}(d_{1}a,d_{2}a,...,d_{n}a)=C^{(n)}\left( d_{1}a\otimes
d_{2}a\otimes ...\otimes d_{n}a\right) ,$

\strut

for all $n\in \mathbb{N},$ and for $d_{1},...,d_{n}\in D_{G},$ where $%
\widehat{C}=(C^{(n)})_{n=1}^{\infty }\in I^{c}\left( W^{*}(G),D_{G}\right) $
is the cumulant multiplicative bimodule map induced by the conditional
expectation $E,$ in the sense of Speicher. We define the $n$-th trivial
moment of $a$ and the $n$-th trivial cumulant of $a$ by

\strut

$\ \ \ \ \ E(a^{n})$ $\ \ $and $\ \ k_{n}\left( \underset{n-times}{%
\underbrace{a,a,...,a}}\right) =C^{(n)}\left( a\otimes a\otimes ...\otimes
a\right) ,$

\strut

respectively, for all $n\in \mathbb{N}.$
\end{definition}

\strut

In [14], we showed that

\strut

\begin{theorem}
(See [14]) Let $n\in \mathbb{N}$ and let $%
L_{w_{1}}^{u_{1}},...,L_{w_{n}}^{u_{n}}\in \left( W^{*}(G),E\right) $ be $%
D_{G}$-valued random variables, where $w_{1},$ $...,$ $w_{n}\in FP(G)$ and $%
u_{j}$ $\in $ $\{1,$ $*\},$ $j$ $=$ $1,$ $...,$ $n.$ Then

\strut

$\ \ \ \ \ \ \ k_{n}\left( L_{w_{1}}^{u_{1}}...L_{w_{n}}^{u_{n}}~\right)
=\mu _{w_{1},...,w_{n}}^{u_{1},...,u_{n}}\cdot
E(L_{w_{1}}^{u_{1}},...,L_{w_{n}}^{u_{n}}),$

\strut

where $\mu _{w_{1},...,w_{n}}^{u_{1},...,u_{n}}=\underset{\pi \in
C_{w_{1},...,w_{n}}^{u_{1},...,u_{n}}}{\sum }\mu (\pi ,1_{n}).$ Here, $%
C_{w_{1},...,w_{n}}^{u_{1},...,u_{n}}$ is a subset of $NC(n)$ consisting of
all partitions $\pi $ in $NC(n),$ satisfying that

\strut

$\ \ \ \ \ \ \ \ \ E_{\pi }\left(
L_{w_{1}}^{u_{1}},...,L_{w_{n}}^{u_{n}}\right)
=E(L_{w_{1}}^{u_{1}}...L_{w_{n}}^{u_{n}})\neq 0_{D_{G}}.$ \ 

$\square $
\end{theorem}

\bigskip \strut \strut \strut \strut

\strut The above theorem show us that, different from the general case, the
mixed $n$-th $D_{G}$-cumulants of operator-valued random variables in $%
(W^{*}(G),$ $E)$ is the product of certain complex number and the mixed $n$%
-th $D_{G}$-moments of the operators. Now, we consider the $D_{G}$-valued
freeness of given two random variables in $\left( W^{*}(G),E\right) $. We
will characterize the $D_{G}$-freeness of $D_{G}$-valued random variables $%
L_{w_{1}}$ and $L_{w_{2}},$ where $w_{1}\neq w_{2}\in FP(G).$ And then we
will observe the $D_{G}$-freeness of arbitrary two $D_{G}$-valued random
variables $a_{1}$ and $a_{2}$ in terms of their supports.

\strut

\begin{definition}
Let $w_{1}$ and $w_{2}$ be elements in the free semigroupoid $\Bbb{F}%
^{+}(G). $ We say that they are diagram-distinct if they have the different
diagrams on $G.$ Also, we will say that the subsets $X_{1}$ and $X_{2}$ of $%
\Bbb{F}^{+}(G)$ are diagram-distinct if $w_{1}$ and $w_{2}$ are
diagram-distinct, for all pair $(w_{1},$ $w_{2})$ in $X_{1}$ $\times $ $%
X_{2}.$
\end{definition}

\strut \strut \strut

By the previous theorem, we can get the following theorem which shows that
the diagram-distinctness characterize the $D_{G}$-freeness of generators of $%
W^{*}(G)$ ;

\strut

\begin{theorem}
(See [14]) Let $w_{1},w_{2}\in FP(G)$ be finite paths. The $D_{G}$-valued
random variables $L_{w_{1}}$ and $L_{w_{2}}$ in $\left( W^{*}(G),E\right) $
are free over $D_{G}$ if and only if $w_{1}$ and $w_{2}$ are
diagram-distinct. $\square $
\end{theorem}

\strut \strut \strut \strut \strut

\begin{corollary}
(See [14]) Let $a,b\in \left( W^{*}(G),E\right) $ be $D_{G}$-valued random
variables with their supports $\mathbb{F}^{+}(G:a)$ and $\mathbb{F}%
^{+}(G:b). $ The $D_{G}$-valued random variables $a$ and $b$ are free over $%
D_{G}$ in $\left( W^{*}(G),E\right) $ if $FP(G:a_{1})$ and $FP(G:a_{2})$ are
diagram-distinct. $\square $
\end{corollary}

\strut \strut \strut \strut

In [15], we observed certain kind of $D_{G}$-valued random variables in $%
\left( W^{*}(G),E\right) .$ One of the most interesting elements are $D_{G}$%
-semicircular elements.

\strut

\begin{proposition}
(See [15]) Let $l\in loop(G)$ be a loop. Then the $D_{G}$-valued random
variable $L_{l}+L_{l}^{*}$ is $D_{G}$-semicircular. $\square $
\end{proposition}

\strut \strut

\strut \strut

\strut

\section{One-Vertex Graph $W^{*}$-Probability Spaces}

\strut

\strut

Throughout this chapter, fix $N\in \Bbb{N}.$ Suppose that $G$ be a finite
directed graph with only one vertex. Let

\strut

\begin{center}
$V(G)=\{v\}$ \ \ and \ \ $E(G)=\{l_{1},...,l_{N}\},$
\end{center}

\strut

where $l_{j}=vl_{j}v$ is a loop, for all \ $j=1,...,N.$ Notice that, in this
case, the diagonal subalgebra $D_{G}$ is isomorphic to $\Bbb{C}.$ i.e,

$\strut $

\begin{center}
$D_{G}=\overline{\Bbb{C}[L_{v}]}^{w}=\Bbb{C}.$
\end{center}

\strut

Thus the canonical conditional expectation $E$ $:$ $W^{*}(G)$ $\rightarrow $ 
$D_{G}$ is a linear map and hence the graph $W^{*}$-probability space $%
\left( W^{*}(G),\text{ }E\right) $ over its diagonal subalgebra $D_{G}$ is
just a (scalar-valued) $W^{*}$-probability space. We will denote such linear
map $E$ by $tr$ and, by $\left( W^{*}(G),\text{ }tr\right) ,$ we will denote
the corresponding graph $W^{*}$-probability space $\left( W^{*}(G),\text{ }%
E\right) .$ Remark that the linear functional $tr=E$ is a faithful trace on $%
W^{*}(G).$ Indeed, assume that $tr\left( xx^{*}\right) $ $=$ $0_{D_{G}}$ $=$ 
$0.$ Then $x$ $=$ $0.$ Also, indeed, $tr$ is a trace. Thus, our graph $W^{*}$%
-probability space $\left( W^{*}(G),tr\right) $ is a tracial $W^{*}$%
-probability space.

\strut \strut

In this setting, the projection $L_{v}$ is the identity $1_{W^{*}(G)}$ of $%
W^{*}(G)$ and the partial isometries $L_{l_{j}},$ \ $j$ $=$ $1,$ $...,$ $N,$
are unitaries in this graph $W^{*}$-algebra, $W^{*}(G),$ since

\strut

\begin{center}
$L_{l_{j}}^{*}L_{l_{j}}=L_{v}=1_{W^{*}(G)}=L_{l_{j}}L_{l_{j}}^{*},$
\end{center}

\strut

for all \ $j=1,...,N,$ and hence

\strut

\begin{center}
$L_{l_{j}}^{*}=L_{l_{j}}^{-1},$ \ for all \ $j=1,...,N.$
\end{center}

\strut

Therefore, the graph $W^{*}$-algebra $W^{*}(G)$ can be understood as a $%
W^{*} $-algebra generated by $N$-unitaries. By [15], we have $D_{G}$%
-semicircular elements $L_{l_{j}}$ $+$ $L_{l_{j}}^{*},$ \ $j$ $=$ $1,$ $...,$
$N.$ Since $D_{G}$ $=$ $\Bbb{C},$ in our case, they are indeed
(scalar-valued) semicircular elements in $\left( W^{*}(G),\text{ }tr\right) .
$ So, we can consider the $W^{*}$-subalgebra $L(G)$ of $W^{*}(G)$ generated
by semicircular elements $\frac{1}{\sqrt{2}}\left( L_{l_{j}}\text{ }+\text{ }%
L_{l_{j}}^{*}\right) $'s, \ $j$ $=$ $1,$ $...,$ $N.$

\strut

\begin{definition}
Let $G$ be the given one-vertex-$N$-loop-edges directed\ graph. Define a $%
W^{*}$-subalgeba $L(G)$ of the graph $W^{*}$-algebra $W^{*}(G)$ by

\strut

$\ \ \ \ \ \ \ \ \ \ L(G)\overset{def}{=}\overline{\Bbb{C}[\{\frac{1}{\sqrt{2%
}}\left( L_{l_{j}}+L_{l_{j}}^{*}\right) :j=1,...,N]}^{w}.$

\strut

Let $\left( W^{*}(G),tr\right) $ be the graph $W^{*}$-probability space
(over $D_{G}$ $=$ $\Bbb{C}$), with its faithful trace $tr$ $=$ $tr\mid
_{L(G)}.$ We will call the $W^{*}$-probability space $\left( L(G),\text{ }%
tr\right) ,$ the semicircular algebra.
\end{definition}

\strut \strut \strut \strut \strut

By [9], we can see that $\left( L(G),tr\right) $ is isomorphic to $\left(
L(F_{N}),\tau \right) .$ Recall that $L(F_{N})$ is the free group factor
induced by the free group $F_{N}$ with $N$-generators. i.e, $%
F_{N}=\,<g_{1},...,g_{N}>$ \ and \ $L(F_{N})=\overline{\Bbb{C}[F_{N}]}^{w}.$
So, if $x\in L(F_{N}),$ then $x$ can be expressed as $x=\underset{g\in F_{N}%
}{\sum }\alpha _{g}g.$ We can define the trace $\tau $ on the free group
factor $L(F_{N})$ by

\strut

\begin{center}
$\tau :L(F_{N})\rightarrow \Bbb{C},$ \ \ $\tau \left( \underset{g\in F_{N}}{%
\sum }\alpha _{g}g\right) =\alpha _{e},$
\end{center}

\strut

where $e$ is the group identity of $F_{N}.$ Voiculescu showed that there
exists a semicircular system $\{x_{j}$ $:$ $j$ $=$ $1,$ $...,$ $N\}$, in the
sense of the set consisting of mutually free semicircular elements with
covariance $1,$ such that

\strut

\begin{center}
$L(F_{N})=vN\left( \{x_{j}:j=1,...,N\}\right) ,$
\end{center}

\strut

where $vN(S)$ means the von Neumann algebra generated by the set $S.$ Define
an operator

\strut

\begin{center}
$T_{0}=\sum_{j=1}^{N}x_{j},$
\end{center}

\strut

where $x_{j}$'s are semicircular elements generating the free group factor $%
L(F_{N}).$ We will call this operator $T_{0}$ the generating operator of $%
L(F_{N}).$

\strut

We can show that the semicircular algebra $\left( L(G),tr\right) $ has the
same free probabilistic structure with $\left( L(F_{N}),\tau \right) .$ i.e,
there exists a $W^{*}$-algebra isomorphism between $L(G)$ and $L(F_{N}),$
which preserves the moments of all generators (See [9]). It is easy to do
that by defining the generator-preserving linear map between $L(G)$ and $%
L(F_{N})$, by regarding $L(F_{N})$ as the von Neumann algebra $vN$ $(\{x_{j}$
$:$ $j$ $=$ $1,$ $...,$ $N\}).$

\strut

\begin{proposition}
Let $G$ be the given one-vertex-$N$-loop-edges directed graph. The
semicircular algebra $\left( L(G),tr\right) ,$ generated by the semicircular
system

\strut \strut

$\ \ \ \ \ \ \ \ \ \ \ \ \ \{\frac{1}{\sqrt{2}}\left(
L_{l_{j}}+L_{l_{j}}^{*}\right) :j=1,...,N\}$

\strut

is isomorphic to $\left( L(F_{N}),\tau \right) ,$ in the sense of [9]. $%
\square $
\end{proposition}

\strut \strut \strut \strut

Now, we will consider the generating operator contained in $\left(
L(G),tr\right) .$

\strut

\begin{definition}
Let $T\in \left( L(G),tr\right) $ be a random variable defined by

\strut

$\ \ \ \ \ \ \ \ \ \ \ \ \ \ T=\sum_{j=1}^{N}\left( \frac{1}{\sqrt{2}}\left(
L_{l_{j}}+L_{l_{j}}^{*}\right) \right) ,$

\strut

where $\{\frac{1}{\sqrt{2}}\left( L_{l_{j}}+L_{l_{j}}^{*}\right)
:j=1,...,N\} $ is the generator set of $L(G).$ We will call $T$ the
generating opeartor of $L(G).$
\end{definition}

\strut

Recall the generating operator $T_{0}$ of the free group factor $L(F_{N}),$ $%
T_{0}=\sum_{j=1}^{N}x_{j}.$ By [9] and by the previous proposition, we have
the following result;

\strut

\begin{proposition}
Let $T_{0}=\sum_{j=1}^{N}x_{j}$ be the generating operator of $\left(
L(F_{N}),\tau \right) $ and let $T=\sum_{j=1}^{N}\left(
L_{l_{j}}+L_{l_{j}}^{*}\right) $ be the generating operator of $\left(
L(G),tr\right) ,$ where $G$ is the given one-vertex-$N$-loop-edge directed
graph. Then

\strut

$\ \ \ \ \ \ \ \ \ \ \ \ \ \ \ \ \ \ \ \ \ \tau (T_{0}^{n})=tr\left(
T^{n}\right) $

and hence

$\ \ \ \ \ \ \ \ \ \ \ \ \ \ \ \ \ k_{n}^{\tau }\left(
T_{0},...,T_{0}\right) =k_{n}\left( T,..,T\right) ,$

\strut

for all $n\in \Bbb{N}$, where $k_{n}^{\tau }(...)$ and $k_{n}(...)$ are
cumulants with respect to the traces $\tau $ and $tr,$ respectively. In
other words, the operators $T_{0}$ and $T$ are identically distributed. $%
\square $
\end{proposition}

\strut

Now, let $\left( A,\varphi \right) $ be a $W^{*}$-probability space and let $%
a\in \left( A,\varphi \right) $ be a random variable. Then we can define the 
$n$-th moments and the $n$-th cumulants of $a$ by

\strut

\begin{center}
$\varphi (a^{n})$ \ \ \ and \ \ \ $k_{n}^{\varphi }\left( a,...,a\right) ,$
\end{center}

\strut

for all $n\in \Bbb{N},$ where $k_{n}^{\varphi }(...)$ is the cumulant
function induced by $\varphi .$ Defin $\Theta _{1}$ as a set of all formal
power series in the indeterminent $z,$ without the constant terms, in $\Bbb{C%
}[[z]],$ where $\Bbb{C}[[z]]$ is the set of all formal power series. For the
given random variable $a\in \left( A,\varphi \right) ,$ we can define the
following two elements in $\Theta _{1}$ ;

\strut

\begin{center}
$M_{a}(z)=\sum_{n=1}^{\infty }\varphi (a^{n})\,z^{n}$
\end{center}

and

\begin{center}
$R_{a}(z)=\sum_{n=1}^{\infty }\,k_{n}^{\varphi }\left( a,...,a\right)
\,z^{n},$
\end{center}

\strut

called the moment series of $a$ and the R-transform of $a,$ respectively. By
the previous proposition, we can get that ;

\strut

\begin{corollary}
Let $T_{0}$ and $T$ be given as before. Then

\strut

$\ \ \ \ \ \ \ \ \ M_{T_{0}}(z)=M_{T}(z)$ \ \ and \ \ $%
R_{T_{0}}(z)=R_{T}(z), $

\strut

in $\Theta _{1}.$ \ $\square $
\end{corollary}

\strut

Again, let $\left( A_{i},\varphi _{i}\right) $ be a $W^{*}$-probability
space, for $i=1,2,$ and let $a_{i}\in \left( A_{i},\varphi _{i}\right) $ be
random variables, for $i=1,2.$ We say that the random variables $a_{1}$ and $%
a_{2}$ are identically distributed if their R-transforms are same in $\Theta
_{1}.$ i.e, the random variables $a_{1}$ and $a_{2}$ are identically
distributed if

\strut

\begin{center}
$R_{a_{1}}(z)=R_{a_{2}}(z)$ \ \ in \ $\Theta _{1}.$
\end{center}

\strut

Notice that, by the M\"{o}bius inversion, if $a_{1}$ and $a_{2}$ are
identically distributed, then

\strut

\begin{center}
$M_{a_{1}}(z)=M_{a_{2}}(z)$ \ \ in \ \ $\Theta _{1}.$
\end{center}

\strut

The above corollary says that, as random variables, the generating operators 
$T_{0}\in \left( L(F_{N}),\tau \right) $ and $T\in \left( L(G),tr\right) $
are identically distributed. Therefore, by computing the moment series or
R-transform of $T_{0},$ we can get those of $T.$ So, by using the graph $%
W^{*}$-probability technique, we can get the moment series and the
R-transform of $T_{0}\in \left( L(F_{N}),\tau \right) .$

\strut

\strut \strut \strut

\strut

\section{Moment and Cumulants of $T_{0}$}

\strut

\strut

\strut

Throughout this chapter, fix $N\in \Bbb{N}$ and let $G$ be a one-vertex-$N$%
-loop-edge directed graph with

\strut

\begin{center}
$V(G)=\{v\}$ \ \ and \ \ $E(G)=\{l_{j}=vl_{j}v:\,j=1,...,N\}.$
\end{center}

\strut

Recall that $L_{v}=1_{W^{*}(G)}=1_{L(G)}$ and $L_{l_{j}}$'s are unitaries in 
$W^{*}(G),$ for all $j=1,...,N.$ Also, let $\left( W^{*}(G),tr\right) $ be
the graph $W^{*}$-probability space (over its diagonal subalgebra $D_{G}=%
\Bbb{C}$), with its faithful trace on $W^{*}(G).$ For the random variables $%
L_{l_{1}}+L_{l_{1}}^{*},$ ..., $L_{l_{N}}+L_{l_{N}}^{*},$ we can form the
semicircular system and then we can construct the semicircular algebra $%
\left( L(G),tr\right) ,$ defined by

$\strut $

\begin{center}
$L(G)\overset{def}{=}\overline{\Bbb{C}[\{\frac{1}{\sqrt{2}}\left(
L_{l_{j}}+L_{l_{j}}^{*}\right) :j=1,...,N\}]}^{w}$
\end{center}

and

\begin{center}
$tr=tr\mid _{L(G)}.$
\end{center}

\strut

Again, remark that $\left( L(G),tr\right) =\left( L(F_{N}),\tau \right) ,$
where $\left( L(F_{N}),\tau \right) $ is the free group factor induced by
the free group $F_{N},$ with $N$-generators. Notice that, by regarding $%
L(F_{N})$ as the von Neumann algebra $vN\left( \{x_{j}:j=1,...,N\}\right) ,$
generated by the semicircular system $\{x_{1},...,x_{N}\},$ we can get the
above equality, by Voiculescu. In this chapter, we will compute the moments
and cumulants of the generating operator

\strut

\begin{center}
$T=\sum_{j=1}^{N}\left( L_{l_{j}}+L_{l_{j}}^{*}\right) $ \ of \ $L(G).$
\end{center}

\strut

Since the generating operator $T_{0}=\sum_{j=1}^{N}x_{j}$ of $\left(
L(F_{N}),\tau \right) $ and the generating operator $T$ of $\left(
L(G),tr\right) $ are identically distributed, the computations for $T$ will
be the reformulation of the moments and cumulants of $T_{0}.$ In fact, the
moments and cumulants of such element $T_{0}$ is solved in various articles.
However, in this section, we will provides the graph $W^{*}$-probability
approach.

\strut \strut

\begin{theorem}
Let $G$ be the given one-vertex-$N$-loop-edge directed graph and let $\left(
L(G),tr\right) $ be the semicircular algebra generated by semicircular
elements $\frac{1}{\sqrt{2}}\left( L_{l_{j}}+L_{l_{j}}^{*}\right) ,$ $%
j=1,...,N.$ If $T=\sum_{j=1}^{N}\frac{1}{\sqrt{2}}\left(
L_{l_{j}}+L_{l_{j}}^{*}\right) $ is the generating operator of $L(G),$ then
it has all vanishing odd moments and cumulants and

\strut

(1) $\ \ tr\left( T^{n}\right) =c_{\frac{n}{2}}\cdot N^{\frac{n}{2}},$

\strut

(2) \ $k_{n}^{tr}\left( T,...,T\right) =\left\{ 
\begin{array}{lll}
N &  & \text{if }n=2 \\ 
&  &  \\ 
0 &  & \text{otherwise,}
\end{array}
\right. ,$

\strut

for all $n\in 2\Bbb{N},$ where $c_{k}=\frac{1}{k+1}\left( 
\begin{array}{l}
2k \\ 
\,\,k
\end{array}
\right) $ is the $k$-th Catalan number.
\end{theorem}

\strut

\begin{proof}
Fix $n\in \Bbb{N}.$ If $n$ is odd, then we have the vanishing moments of $T,$
because of the $*$-axis-property. Thus all odd cumulants of $T$ also vanish.
We will prove (2), first.

\strut

(2) Notice that $L_{l_{1}}+L_{l_{1}}^{*},$ ..., $L_{l_{N}}+L_{l_{N}}^{*}$
are free from each other in $\left( L(G),tr\right) ,$ by the
diagram-dsistinctness of $l_{1},$ ..., $l_{N}$. So, we have that

\strut

$\ k_{n}^{tr}\left( T,...,T\right) =k_{n}^{tr}\left( \sum_{j=1}^{N}\frac{1}{%
\sqrt{2}}(L_{l_{j}}+L_{l_{j}}^{*}),\,...,\,\sum_{j=1}^{N}\frac{1}{\sqrt{2}}%
(L_{l_{j}}+L_{l_{j}}^{*})\right) $

\strut

$\ \ \ =\sum_{j=1}^{N}k_{n}^{tr}\left( \frac{1}{\sqrt{2}}\left(
L_{l_{j}}+L_{l_{j}}^{*}\right) ,...,\frac{1}{\sqrt{2}}\left(
L_{l_{j}}+L_{l_{j}}^{*}\right) \right) $

\strut

by the mutually freeness of $L_{l_{1}}+L_{l_{1}}^{*},$ ..., $%
L_{l_{N}}+L_{l_{N}}^{*}$

\strut

$\ \ \ =\sum_{j=1}^{N}\,\underset{(u_{1},...,u_{n})\in \{1,*\}^{n}}{\sum }%
k_{n}^{tr}\left( \frac{1}{\sqrt{2}}L_{l_{j}}^{u_{1}},...,\frac{1}{\sqrt{2}}%
L_{l_{j}}^{u_{n}}\right) $

\strut

$\ \ \ =\left\{ 
\begin{array}{ll}
\sum_{j=1}^{N}\,\underset{(u_{1},u_{2})\in \{1,*\}^{2}}{\sum }%
k_{n}^{tr}\left( \frac{1}{\sqrt{2}}L_{l_{j}}^{u_{1}},\frac{1}{\sqrt{2}}%
\,L_{l_{j}}^{u_{2}}\right) & \text{if }n=2 \\ 
&  \\ 
0_{D_{G}}=0 & \text{otherwise}
\end{array}
\right. $

\strut

by the semicircularity of $L_{l_{j}}+L_{l_{j}}^{*},$ for all \ $j=1,...,N$

\strut

$\ \ \ =\left\{ 
\begin{array}{ll}
\sum_{j=1}^{N}\,\frac{1}{2}\left( 2L_{v}\right) & \text{if }n=2 \\ 
&  \\ 
0 & \text{otherwise}
\end{array}
\right. $

\strut

by Section 2.5

\strut

$\ \ \ =\left\{ 
\begin{array}{ll}
\sum_{j=1}^{N}\,1=N & \text{if }n=2 \\ 
&  \\ 
0 & \text{otherwise,}
\end{array}
\right. $

\strut

since $L_{v}=1_{L(G)}=1\in \Bbb{C}.$

\strut

(1) Now, remark that the generating operator $T$ is semicircular, by (2).
Fix $n\in 2\Bbb{N}.$ Then we have that

\strut

$\ \ \ \ \ tr\left( T^{n}\right) =\underset{\pi \in NC(n)}{\sum }k_{\pi
}\left( T,...,T\right) $

\strut

by the M\"{o}bius inversion

\strut

$\ \ \ \ \ \ \ \ \ \ \ \ \ =\underset{\pi \in NC_{2}(n)}{\sum }k_{\pi
}(T,...,T)$

\strut

by the semicircularity of $T,$ where

\strut

$\ \ \ \ \ \ \ NC_{2}(n)=\{\pi \in NC(n):V\in \pi \Rightarrow \left|
V\right| =2\},$

\strut

and then

\strut

$\ \ \ \ \ \ \ \ \ \ \ \ =\underset{\pi \in NC_{2}(n)}{\sum }\left( 
\underset{V\in \pi }{\prod }k_{V}(T,...,T)\right) $

\strut

$\ \ \ \ \ \ \ \ \ \ \ \ =\underset{\pi \in NC_{2}(n)}{\sum }\left(
k_{2}(T,T)\right) ^{\left| \pi \right| }=\underset{\pi \in NC_{2}(n)}{\sum }%
N^{\left| \pi \right| }$

\strut

since $k_{2}(T,T)=N,$ by (2)

\strut

$\ \ \ \ \ \ \ \ \ \ \ =\underset{\pi \in NC_{2}(n)}{\sum }N^{\frac{n}{2}%
}=\left| NC_{2}(n)\right| \cdot N^{\frac{n}{2}}=c_{\frac{n}{2}}\cdot N^{%
\frac{n}{2}},$

\strut

since $\left| NC_{2}(n)\right| =\left| NC(\frac{n}{2})\right| =c_{\frac{n}{2}%
},$ where $c_{k}$ is the $k$-th Catalan number.
\end{proof}

\strut

By the previous theorem we can get that ;

\strut

\begin{corollary}
Let $T$ be the generating operator of the semicircular algebra $\left(
L(G),tr\right) ,$ where $G$ is the given one-vertex-$N$-loop-edge directed
graph. Then the moment sereis $M_{T}(z)$ of $T$ and the R-transform $%
R_{T}(z) $ of $T$ are

\strut

$\ \ \ \ \ \ \ \ \ \ \ \ \ \ \ \ M_{T}(z)=$ $\sum_{n=1}^{\infty }\left( c_{%
\frac{n}{2}}\cdot N^{\frac{n}{2}}\right) z^{n}$

and

$\ \ \ \ \ \ \ \ \ \ \ \ \ \ \ \ \ \ \ \ \ \ \ \ R_{T}(z)=N\cdot z^{2},$

\strut

in $\Theta _{1}.$ \ $\square $
\end{corollary}

\strut \strut

The above corollary shows that the generating operator $T_{0}$ of the free
group factor $L(F_{N})$ satisfies that

\strut

\begin{center}
$M_{T_{0}}(z)=$ $\sum_{n=1}^{\infty }\left( c_{\frac{n}{2}}\cdot N^{\frac{n}{%
2}}\right) z^{n}$
\end{center}

and

\begin{center}
$R_{T_{0}}(z)=N\cdot z^{2},$
\end{center}

in $\Bbb{C}[[z]],$ too.

\strut \strut

\strut

\strut

\section{Embedding $L(F_{N})$ into $W^{*}(G)$}

\strut

\strut

\strut

In this chapter, we will consider the embedding of the free group factor $%
L(F_{N})$ in the graph $W^{*}$-probability space $\left( W^{*}(G),E\right) ,$
where $G$ is an arbitrary countable directed graph having at least one
vertex with $N$-diagram-distinct loops. This is already observed in [14].
Throughout this chapter, let $G$ be a countable directed graph and let $%
\left( W^{*}(G),E\right) $ be the graph $W^{*}$-probability space over its
diagonal subalgebra $D_{G}.$ Also, we will assume that there exists a vertex 
$v_{0}\in V(G)$ such that there is a nonempty set consisting of basic loops,

\strut

\begin{center}
$Loop_{v_{0}}(G)=\{l\in Loop(G):l=v_{0}lv_{0}\}$
\end{center}

\strut and

\begin{center}
$\left| Loop_{v_{0}}(G)\right| =N.$
\end{center}

\strut

Without loss of generality, we can write $Loop_{v_{0}}(G)$ $=$ $\{l_{1},$ $%
...,$ $l_{N}\}.$ Notice that $l_{1},$ $...,$ $l_{N}$ are distinct basic
loops. So, we can construct the $D_{G}$-semicircular system,

\strut

\begin{center}
$S\overset{def}{=}\{\frac{1}{\sqrt{2}}\left( L_{l}+L_{l}^{*}\right) :l\in
Loop_{v_{0}}(G)\}.$
\end{center}

\strut

i.e, the set $S$ is consist of mutually $D_{G}$-free $D_{G}$-semicircular
elements in $\left( W^{*}(G),E\right) .$ Now, define the (scalar-valued)
semicircular subalgebra $L(S)$ by

\strut

\begin{center}
$L(S)=\overline{\Bbb{C}[S]}^{w}.$
\end{center}

\strut

Notice that this subalgebra $L(S)$ is slightly different from those of [16].
In [16], we defined the $D_{G}$-semicircular subalgebra $L_{D_{G}}(S)$ by

\strut

\begin{center}
$L_{D_{G}}\left( S\right) =\overline{D_{G}[S]}^{w}.$
\end{center}

\strut

We can see that

\strut

\begin{center}
$\left( L_{D_{G}}(S),E\right) =\left( L(S),E\mid _{L(S)}\right) \otimes
\left( D_{G},\mathbf{1}\right) ,$
\end{center}

\strut

where $\mathbf{1}$ is the identity map on $D_{G}.$ More generally, we have
that;

\strut 

\begin{theorem}
Let 

\strut 

$\ \mathcal{L}_{N}=\left\{ \frac{1}{\sqrt{2}}\left(
L_{l_{j}}+L_{l_{j}}^{*}\right) :
\begin{array}{l}
l_{j}=v_{j}l_{j}v_{j},j=1,...,N \\ 
l_{j}\text{'s are mutually diagram-distinct}
\end{array}
\right\} $

\strut 

be a $D_{G}$-valued semicircular system in $W^{*}(G).$ Then

\strut 

$\ \ \ \ \ \ \left( vN(\mathcal{L}_{N},D_{G}),\text{ }E\right) =\left(
vN\left( \mathcal{L}_{N},D_{N}\right) ,E\right) \otimes \left( D_{G},\mathbf{%
1}\right) ,$

\strut \strut 

where $vN(S)$ is the von Neumann algebra generated by the set $S$ and $%
\mathbf{1}$ is the identity map on $D_{G},$ and 

\strut 

\ \  $\ \ \ \ \ \ \ D_{N}=\overline{\Bbb{C}%
[\{L_{v_{j}}:l_{j}=v_{j}l_{j}v_{j},\text{ }j=1,...,N\}]}^{w}.$
\end{theorem}

\strut \strut \strut 

\begin{proof}
Let $\mathcal{L}_{N}$ be the collection of $N$-$D_{G}$-semicircular elements
which are mutually free over $D_{G},$ as follows;

\strut 

$\ \ \ \mathcal{L}_{N}=\left\{ \frac{1}{\sqrt{2}}\left(
L_{l_{j}}+L_{l_{j}}^{*}\right) :
\begin{array}{l}
l_{j}=v_{j}l_{j}v_{j},j=1,...,N \\ 
l_{j}\text{'s are mutually diagram-distinct}
\end{array}
\right\} ,$

\strut 

The $\mathcal{L}_{N}$ is a $D_{G}$-semicircular system. As $W^{*}$-algebras, 

\strut 

$\ \ \ \ \ \ \ \ \ \ \ \ \ vN$ $(\mathcal{L}_{N},$ $D_{G})$ $\simeq $ $%
vN\left( \mathcal{L}_{N},\text{ }D_{N}\right) $ $\otimes $ $D_{G},$

\strut 

where $D_{N}=\overline{C[\{L_{v_{j}}:l_{j}=v_{j}l_{j}v_{j}\}]}^{w}.$ Indeed,
without loss of generality, take $a$ $\in $ $vN(\mathcal{L}_{N},$ $D_{G})$ by

$\strut $

$\ \ \ \ \ \
a=d_{1}a_{l_{i_{1}}}^{k_{1}}d_{2}a_{l_{i_{2}}}^{k_{2}}...d_{n}a_{l_{i_{n}}}^{k_{n}}
$ \ and \ $a_{l_{j}}=\frac{1}{\sqrt{2}}\left( L_{l_{j}}+L_{l_{j}}^{*}\right) 
$

\strut

where $d_{1},...,d_{n}\in D_{G},$ $k_{1},...,k_{n}\in \Bbb{N}$ and $(i_{1},$ 
$...,$ $i_{n})$ $\in $ $\{1,$ $...,$ $N\}^{n},$ $n$ $\in $ $\mathbb{N}.$
Observe that, for any $j$ $\in $ $\{1,$ $...,$ $N\},$ we have that

$\strut $

$\ \ \ \ \ \ \ \ \ a_{l_{j}}^{k}=\left( L_{l_{j}}+L_{l_{j}}^{*}\right)
^{k}=L_{l_{j}^{k}}+L_{l_{j}^{k}}^{*}+Q\left( L_{l_{j}},L_{l_{j}}^{*}\right)
, $

\strut

where $Q\in \mathbb{C}[z_{1},z_{2}]$. Also, observe that $%
L_{l_{j}}^{k_{1}}L_{l_{j}}^{*\,\,k_{2}},$ for any $k_{1}$ $,k_{2}$ $\in $ $%
\mathbb{N},$ satisfies that

\strut

$\ \ \
L_{l_{j}}^{k_{1}}L_{l_{j}}^{*\,%
\,k_{2}}=L_{l_{j}^{k_{1}}}L_{l_{j}^{k_{2}}}^{*}=\left\{ 
\begin{array}{ll}
L_{l_{j}^{k_{1}-k_{2}}}=L_{v_{j}}L_{l_{j}^{k_{1}-k_{2}}} & \text{if }%
k_{1}>k_{2} \\ 
L_{l_{j}^{k_{2}-k_{1}}}^{*}=L_{v_{j}}L_{l_{j}^{k_{2}-k_{1}}}^{*} & \text{if }%
k_{1}<k_{2} \\ 
L_{v_{j}} & \text{if }k_{1}=k_{2},
\end{array}
\right. $

and similarly,

\strut

$\ \ \
L_{l_{j}}^{*\,\,k_{1}}L_{l_{j}}^{%
\,k_{2}}=L_{l_{j}^{k_{1}}}^{*}L_{l_{j}^{k_{2}}}=\left\{ 
\begin{array}{ll}
L_{l_{j}^{k_{1}-k_{2}}}^{*}=L_{v_{j}}L_{l_{j}^{k_{1}-k_{2}}}^{*} & \text{if }%
k_{1}>k_{2} \\ 
L_{l_{j}^{k_{2}-k_{1}}}=L_{v_{j}}L_{l_{j}^{k_{2}-k_{1}}} & \text{if }%
k_{1}<k_{2} \\ 
L_{v_{j}} & \text{if }k_{1}=k_{2}.
\end{array}
\right. $

So,

\strut

$\ \ \ \ \ \ \ \ \ \ \ \ \ \ \ Q(L_{l_{j}},L_{l_{j}}^{*})=L_{v_{j}}\left(
Q(L_{l_{j}},L_{l_{j}}^{*})\right) L_{v_{j}},$

\strut

for all \ $j=1,...,N.$ Thus

\strut

(1.1) $\ \ \ \ \ \ \
a_{l_{j}}^{k}=L_{v_{j}}a_{l_{j}}^{k}=L_{v_{j}}a_{l_{j}}^{k}L_{v_{j}},$ for
all $j=1,...,N.$

\strut \strut \strut

Now, consider that

\strut

$\ \ \ \ \ \ \ \ \ \ \ \ \ d_{j}=d_{j}^{N}+d_{j}^{\prime },$ $\forall
j=1,...,N.$

\strut

where $d_{j}^{N}=\sum_{j=1}^{N}L_{v_{j}}d_{j}L_{v_{j}}$ and $d_{j}^{\prime
}=d_{j}-d_{j}^{N}$ in $D_{G}.$ So, we can rewrite that

\strut

$\ \ a=\left( d_{1}^{N}+d_{1}^{\prime }\right) a_{l_{i_{1}}}^{k_{1}}\left(
d_{2}^{N}+d_{2}^{\prime }\right) a_{l_{i_{2}}}^{k_{2}}...\left(
d_{n}^{N}+d_{n}^{\prime }\right) a_{l_{i_{n}}}^{k_{n}}$

\strut

$\ \ \ \
=d_{1}^{N}a_{l_{i_{1}}}^{k_{1}}d_{2}^{N}a_{l_{i_{2}}}^{k_{2}}...d_{n}^{N}a_{l_{i_{n}}}^{k_{n}}+d_{1}^{\prime }a_{l_{i_{1}}}^{k_{1}}d_{2}^{\prime }a_{l_{i_{2}}}^{k_{2}}...d_{n}^{\prime }a_{l_{i_{n}}}^{k_{n}} 
$

\strut

$\ \ \ \
=d_{1}^{N}a_{l_{i_{1}}}^{k_{1}}d_{2}^{N}a_{l_{i_{2}}}^{k_{2}}...d_{n}^{N}a_{l_{i_{n}}}^{k_{n}}, 
$

\strut \strut

by (1.1). This shows that $a=a\otimes 1\in vN(\mathcal{L}_{N},D_{N})\otimes 1
$ and

\strut

$\ \ \ \ \ \ \ \ E\left( a\right) =E_{D_{N}}^{D_{G}}\circ
E(a)=E_{N}(a)=E_{N}\otimes \mathbf{1}(a\otimes 1).$

\strut

Trivially, if $a\in D_{G}\subset L_{D_{G}}(S),$ then $a$ $=$ $1$ $\otimes $ $%
a$ $\in $ $1$ $\otimes $ $D_{G}.$ Futhermore, if $a$ $\in $ $D_{G},$ then

\strut

\ $\ \ \ \ \ \ \ \ \ \ E(a)=a=1\otimes a=E_{N}\otimes \mathbf{1}(1\otimes
a). $

\strut \strut \strut 
\end{proof}

\strut 

By the previous theorem, as a corollary, we can get that;

\strut 

\begin{corollary}
$\left( L_{D_{G}}(S),E\right) =\left( L(S),E\mid _{L(S)}\right) \otimes
\left( D_{G},\mathbf{1}\right) .$ $\square $
\end{corollary}

\strut \strut \strut 

Therefore, we have that;

\strut

\begin{corollary}
Let $v_{0},$ $S$ and $L(S)$ be given as above. Then $\left( L(S),E\mid
_{L(S)}\right) =\left( L(F_{N}),\tau \right) .$
\end{corollary}

\strut

\begin{proof}
Denote $\frac{1}{\sqrt{2}}\left( L_{l_{j}}+L_{l_{j}}^{*}\right) $ by $x_{j},$
for all $j=1,...,N.$ Then, for any projections $L_{v}$ $\in $ $W^{*}(G),$ $v$
$\in $ $V(G),$ we have that

\strut

$\ \ \ \ \ \ \ \ \ \ \ \ \ L_{v}x_{j}=x_{j}L_{v},$ \ for all \ $j=1,...,N.$

\strut

Suppose that $v\neq v_{0}$ in $V(G).$ Then

\strut

$\ \ \ \ \ \ \ \ \ L_{v}x_{j}=0_{D_{G}}=x_{j}L_{v},$ for all \ $j=1,...,N$

\strut

Now, assume that $v=v_{0}$ in $V(G).$ Then

\strut

$\ \ \ \ \ \ \ \ \ L_{v_{0}}x_{j}=x_{j}=x_{j}L_{v_{0}},$ \ for all \ $%
j=1,...,N.$

\strut

So, the conditional expectation $E$ on $L(S)$ is regarded as the linear
functional $E$ $:$ $L(S)$ $\rightarrow $ $\Bbb{C\xi }_{v_{0}}.$ Moreover
this linear functional $E$ is faithful and tracial. So, the $D_{G}$%
-semicircular elements $x_{j}$'s are semicircular in the $W^{*}$-subalgebra $%
\left( L(S),E\right) .$
\end{proof}

\strut

The above corollary shows how to embed the free group factor $L(F_{N})$ into 
$W^{*}(G).$ Vice versa, if a graph $G$ contains a vertex having $N$-loops
based on it, then we can construct a $W^{*}$-subalgebra $L(S)$ isomorphic to
the free group factor $L(F_{N}).$

\strut

\strut

\strut

\strut \textbf{References}

\strut

\strut

{\small [1] \ \ A. Nica, R-transform in Free Probability, IHP course note.}

{\small [2] \ \ A. Nica, R-transforms of Free Joint Distributions and
Non-crossing Partitions, J. of Func. Anal, 135 (1996), 271-296.\strut }

{\small [3] \ \ A. Nica, D. Shlyakhtenko, R. Speicher, R-Cyclic Families of
Matrices in Free Probability, J. of Funct. Anal, 188 (2002), 227-271.}

{\small [4] \ \ A. Nica, D. Shlyakhtenko, R. Speicher, R-Diagonal Elements
and Freeness with Amalgamation, Canad. J. Math, 53, \# 2, (2001), 335-381.}

{\small [5] \ \ A. Nica, R. Speicher, R-diagonal Pair-A Common Approach to
Haar Unitaries and Circular Elements, (1995), Preprint.}

{\small [6] \ \ D. Shlyakhtenko, Some Applications of Freeness with
Amalgamation, J. Reine Angew. Math, 500 (1998), 191-212.\strut }

{\small [7] \ \ D. Shlyakhtenko, A-Valued Semicircular Systems, J. of Funct
Anal, 166 (1999), 1-47.\strut }

{\small [8] \ \ D. Voiculescu, Operations on Certain Non-commuting
Operator-Valued Random Variables, Ast\'{e}risque, 232 (1995), 243-275.\strut 
}

{\small [9] \ \ D.Voiculescu, K. Dykemma and A. Nica, Free Random Variables,
CRM Monograph Series Vol 1 (1992).\strut }

{\small [10] F. Radulescu, Singularity of the Radial Subalgebra of }$%
L(F_{N}) ${\small \ and the Puk\'{a}nszky Invariant, Pacific J. of Math,
vol. 151, No 2 (1991)\strut , 297-306.\strut \strut }

{\small [11] I. Cho, Amalgamated Boxed Convolution and Amalgamated
R-transform Theory (2002), Preprint.}

{\small [12]\strut I. Cho, An Example of Moment Series under the
Compatibility (2003), Preprint.}

{\small [13] I. Cho, The Moment Series and R-transform of the Generating
Operator of }$L(F_{N})${\small \ (2003), Preprint.}

{\small [14] I. Cho, Graph }$W^{*}${\small -Probability Theory (2004),
Preprint.}

{\small [15] I. Cho, Random Variables in Graph }$W^{*}${\small -Probability
Spaces (2004), Ph.D Thesis, Univ. of Iowa.}

{\small [16] I. Cho, Amalgamated Semicircular Systems in Graph }$W^{*}$%
{\small -Probability Spaces (2004), Preprint.}

{\small [17] I. Cho, Free Product Structure of Graph }$W^{*}${\small %
-Probability Spaces (2004), Preprint.}

{\small [18] R. Speicher, Combinatorial Theory of the Free Product with
Amalgamation and Operator-Valued Free Probability Theory, AMS Mem, Vol 132 ,
Num 627 , (1998).}

{\small [19] R. Speicher, Combinatorics of Free Probability Theory IHP
course note.}

\end{document}